\documentclass{article}
\usepackage{amssymb}
\usepackage{amsmath}

\input{tcilatex}
\begin{document}

\begin{center}
\textbf{A NOTE ON\ INEXTENSIBLE\ FLOWS OF\ PARTIALLY\ AND\ PSEUDO\ NULL\
CURVES\ IN }$E_{1}^{4}$

$^{1,2}$Z\"{u}hal K\"{U}\c{C}\"{U}KARSLAN Y\"{U}ZBA\c{S}I, $^{2}$Mehmet BEKTA%
\c{S}

$^{1}$Brock University, Department of Mathematics, St. Catharines, Ontario
L2S 3A1, Canada

$^{2}$Firat University, Faculty of Science, Department of Mathematics, 23119
Elaz\i \u{g} / Turkey

zkucukarslan@brocku.ca, mbektas@firat.edu.tr
\end{center}

\textbf{Abstract.} In this paper, we study inextensible flows of partially
null and pseudo null curves in $E_{1}^{4}.$ We give neccessary and sufficent
conditions for inextensible flows of partially null and pseudo null curves
in $E_{1}^{4}$



\bigskip \textbf{1. Introduction}

Recently, the study of the motion of inextensible curves has arisen in a
number of diverse engineering applications. The flow of a curve is said to
be inextensible if the arc length is preserved. Physically, inextensible
curve flows give rise to motions in which no strain energy is induced. \ The
swinging motion of a cord of fixed length, for example, or of a piece of
paper carried by the wind, can be described \ by inextensible curve and
surface flows. Such motions arise quite naturally in a wide range of a
physical applications. For example, both Chirikjian and Burdick [6] and
Mochiyama et al. [19] study the shape control of hyper-redundant, or
snake-like robots. Inextensible curve and surface flows also arise in the
context of many problems in computer vision [11], [18] and computer
animation [7], and even structural mechanics [21].

Firstly, Kwon and Park studied inextensible flows of curves and developable
surfaces, which its arclength is preserved, in Euclidean 3-space [16\ ].
Inextensible flows of curves are studied in many different spaces. G\"{u}rb%
\"{u}z have examined inextensible flows of spacelike, timelike and null
curves in [9]. After this work \"{O}grenmi\c{s} et al. have studied
inextensible curves in Galilean space [20] and Y\i ld\i z et al. have
studied inextensible flows of curves according to Darboux frame in Euclidean
3-space [22] . Moreover Latifi et al. (2008) studied inextensible flows of
curves in Minkowski 3-space [17].

In [4], [5], [12] and [13] the authors focused on timelike and space like
curves in $E_{1}^{3}$ and $E_{1}^{4}$ .In the recent work [23], [24] \"{O}%
.G. Y\i ld\i z et al. gave necessary and sufficient conditions for
inextensible flows of non-null curves in $E^{n}$ and $E_{1}^{n}$.

More generally, from the differential geometric point of view, the study of
null curves has its own geometric interest. Many of the classical results
from Riemannian geometry have Lorentz counterparts. In fact, spacelike
curves or timelike curves can be studied by a similar approach to that in
positive definite Riemannian geometry. However, null curves have very
diferent properties from spacelike or timelike curves. In other words, null
curve theory has many results which have no Riemannian analogues. The
presence of null curves often causes important and interesting differences,
as will be the case in the present study [2].

Nowadays, many important and intensive studies are seen about null curves in
Minkowski space. Papers in [1], [3],[8],[14], [15] \ show that how important
field of interest null curves have and obtained some new characterization\
of this curve in Minkowski space.

In the present paper following [13], [16],[ 24], we define inextensible
flows of partially null and pseudo curves in $E_{1}^{4}.$ We give neccessary
and sufficent conditions for inextensible flows of partially null and pseudo
null curves in $E_{1}^{4}.$\qquad 

\qquad \textbf{2. Preliminaries}

Let $E_{1}^{4\text{ }}$denote the 4-dimensional Minkoski space -time i.e.
the Euclidean space $E^{4\text{ \ }}$ with the standart flat metric given by%
\begin{equation*}
\text{ }<~,\text{ }>-\text{ }dx_{1}^{2}\text{~}%
+dx_{2}^{2}+dx_{3}^{2}+dx_{4}^{2}\text{,}
\end{equation*}%
where $\ (x_{1},x_{2},x_{3},x_{4})$ is a rectangular coordinate system of $%
E_{1}^{4\text{ }}.$Since $g$ is an indefinite metric, recall that a vector $v
$ in $E_{1}^{4\text{ }}$can have one of three casual characters: it can be
spacelike if $<v,v>\ >0$ or $v=0,~$timelike if $<v,v><0$ and null
(lightlike) if $<v,v>\ =0$ and $v\neq 0.$ The norm of a vector $v$ is given
by $\parallel v\parallel =\sqrt{\mid <v,v>\mid }$. Therefore, $v$ is unit
vector if $<v,v>=\mp 1.$ Next, vectors $v$ and $w$ are said to orthogonal if 
$<v,w>=0.$ Similarly an arbitrary curve $\ \ \gamma (s)$ can be locally
spacelike,$~$timelike or null (lightlike), if all of its velocity $\ \gamma
^{^{\prime }}(s)$ are respectively spacelike,$~$timelike or null
(lightlike). Next $\ \ \ \gamma (s)~$is a unit speed curve if $\ <\gamma
^{^{\prime }}(s),\ \gamma ^{^{\prime }}(s)>=\mp 1.$

Recall that a spacelike curve in $E_{1}^{4\text{ }}$is called pseudo null
curve or partially null curve, if its principal normal vector is null or its
first binormal is null, respectively [10] . A null curve $\gamma $ is
parametrized by arclength function s, if $<\gamma ^{^{\prime \prime }}(s),\
\gamma ^{^{\prime \prime }}(s)>=1.$ In particular, pseudo null curve or
partially null curve $\gamma (s)~$has unit speed, if $<\gamma ^{^{\prime
}}(s),\ \gamma ^{^{\prime }}(s)>=1.$

In the following we use the notations and concepts from~[10], unless
otherwise stated.

Let $\{T,N,B_{1},B_{2}\}$ be the moving Frenet frame along a curve $\gamma $
in $E_{1}^{4}$, consisting of the tangent, the principal normal, the first
binormal and the second binormal vector fields. Depending on the causal
character of $\gamma $, the Frenet equations have the following forms.

\textbf{Case (a)}. \ If $\gamma $ is partially null curve, the Frenet
formulas gave as ([10]):\bigskip 
\begin{equation}
\left[ 
\begin{array}{c}
T^{^{\prime }} \\ 
N^{^{\prime }} \\ 
B_{1}^{^{\prime }} \\ 
B_{2}^{^{\prime }}%
\end{array}%
\right] =\left[ 
\begin{array}{cccc}
0 & k_{1} & 0 & 0 \\ 
-k_{1} & 0 & k_{2} & 0 \\ 
0 & 0 & k_{3} & 0 \\ 
0 & -k_{2} & 0 & -k_{3}%
\end{array}%
\right] \left[ 
\begin{array}{c}
T \\ 
N \\ 
B_{1} \\ 
B_{2}%
\end{array}%
\right]   \tag{2.1}
\end{equation}%
where the third curvature $k_{3}(s)=0$ for each $s$. Such curve has two
curvatures\ $k_{1}(s)$ and $k_{2}(s)$ and lies fully in a lightlike
hyperplane of $E_{1}^{4\text{ }}.$ In particular, the following equations
hold%
\begin{equation}
<T,T>=<N,N>=1,<B_{1},B_{1}>=<B_{2},B_{2}>=0,  \tag{2.2}
\end{equation}%
\begin{equation}
<T,N>=<T,B_{1}>=<T,B_{2}>=<N,B_{1}>=<N,B_{2}>=0,<B_{1},B_{2}>=1.  \tag{2.3}
\end{equation}

\textbf{Case (b)}. If $\gamma $ is pseudo null curve, the Frenet formulas
are ([10]):%
\begin{equation}
\left[ 
\begin{array}{c}
T^{^{\prime }} \\ 
N^{^{\prime }} \\ 
B_{1}^{^{\prime }} \\ 
B_{2}^{^{\prime }}%
\end{array}%
\right] =\left[ 
\begin{array}{cccc}
0 & k_{1} & 0 & 0 \\ 
0 & 0 & k_{2} & 0 \\ 
0 & k_{3} & 0 & -k_{2} \\ 
-k_{1} & 0 & -k_{3} & 0%
\end{array}%
\right] \left[ 
\begin{array}{c}
T \\ 
N \\ 
B_{1} \\ 
B_{2}%
\end{array}%
\right] ,  \tag{2.4}
\end{equation}%
where the first curvature $k_{1}(s)=0$, if $\gamma $ is straight line, or $%
k_{1}(s)=1$ in all other cases. Such curve has two curvatures $k_{2}(s)$ and 
$k_{3}(s)$ and the following conditions are satisfied:%
\begin{equation}
<T,T>=<B_{1},B_{1}>=1,<N,N>=<B_{2},B_{2}>=0,  \tag{2.5}
\end{equation}%
\begin{equation}
<T,N>=<T,B_{1}>=<T,B_{2}>=<N,B_{1}>=<B_{1},B_{2}>=0,<N,B_{2}>=1.  \tag{2.6}
\end{equation}

\textbf{3. Inextensible Flows of partially null curve in} $E_{1}^{4\text{ }}$

Unless otherwise stated we assume that%
\begin{equation*}
\gamma :[0,l]\times \lbrack 0,w)\rightarrow E_{1}^{4}
\end{equation*}%
is a one parameter family of smooth partially null or pseudo null \ curves
in $E_{1}^{4\text{ }}$, where $l$ is the arclength of the initial curve.
Suppose that $u$ is the curve parametrization variable , $0\leq u\leq l$. If
the speed partially null or pseudo null curves $\gamma $ is given by $%
v=\left\Vert \frac{\partial \gamma }{\partial u}\right\Vert $, then the
arclength of $\gamma $ is given as a function of $u$ by%
\begin{equation}
s(u)=\overset{u}{\underset{0}{\int }}~\left\Vert \frac{\partial \gamma }{%
\partial u}\right\Vert du=\overset{u}{\underset{0}{\int }}~v~du,  \tag{3.1}
\end{equation}%
where

\begin{equation*}
~\left\Vert \frac{\partial \gamma }{\partial u}\right\Vert =\sqrt{\
\left\vert <\frac{\partial \gamma }{\partial u},\frac{\partial \gamma }{%
\partial u}>\right\vert }.
\end{equation*}%
The operator $\frac{\partial }{\partial s}$ is given by%
\begin{equation}
\frac{\partial }{\partial s}=\frac{1}{v}\frac{\partial }{\partial u}, 
\tag{3.2}
\end{equation}%
where $v=\left\Vert \frac{\partial \gamma }{\partial u}\right\Vert .$

In this case; the arclength is as follows $ds=vdu$.

\textbf{Definition 3.1.} Let \ $\gamma $ be a partially null or pseudo null
curves in $E_{1}^{4\text{ }}$ and $\{T,N,B_{1},B_{2}\}$be the Frenet frame
of $\gamma $ in Minkowski space-time. Any flow of the partially null or
pseudo null curves can be expressed as follows%
\begin{equation}
\frac{\partial \gamma }{\partial t}=\beta _{1}T+\beta _{2}N+\beta
_{3}B_{1}+\beta _{4}B_{2}.  \tag{3.3}
\end{equation}%
where, $\beta _{i}$ is the $i^{th}$ scalar speed of the partially null
curves $\gamma $.

Let the arclength variation be

\begin{equation*}
s(u,t)=\overset{u}{\underset{0}{\int }}~v~du.
\end{equation*}

In $E_{1}^{4\text{ }}$ ,the requirement that the partially null or pseudo
null curves not be subject to any elongation or compression can be expressed
by the condition

\begin{equation*}
\frac{\partial }{\partial t}s(u,t)=\overset{u}{\underset{0}{\int }}~\frac{%
\partial v}{\partial t}~du=0,
\end{equation*}%
where $u\in \lbrack 0,l].$

\textbf{Definition 3.2}. Let $\gamma $ be a partially null or pseudo null
curves in $E_{1}^{4\text{ }}$. A partially null curve evolution $\gamma
(u,t) $ and its flow $\frac{\partial \gamma }{\partial t}$ are said to be
inextensible if

\begin{equation}
\frac{\partial }{\partial t}\left\Vert \frac{\partial \gamma }{\partial u}%
\right\Vert =0.  \tag{3.4}
\end{equation}

Before deriving the necessary and sufficient condition for inelastic
partially null or pseudo null curves flow, we need the following lemma.

\textbf{Lemma 3.3.} Let \ $\{T,N,B_{1},B_{2}\}$ be the Frenet frame of a
partially null curve $\gamma $ and 
\begin{equation*}
\frac{\partial \gamma }{\partial t}=\beta _{1}T+\beta _{2}N+\beta
_{3}B_{1}+\beta _{4}B_{2}
\end{equation*}%
be a smooth flow of a partially null curves $\gamma $ in $E_{1}^{4\text{ }}.$
If $\gamma (s)$ be unit speed partially null curve in $E_{1}^{4\text{ }},$
with the curvature $k_{3}(s)=0,$then we have the following equality%
\begin{equation}
\frac{\partial v}{\partial t}=\left( \frac{\partial \beta _{1}}{\partial u}%
-\beta _{2}k_{1}v\right) .\ \ \   \tag{3.5}
\end{equation}%
$\ \ $

\textbf{Proof. }Suppose that $\dfrac{\partial \gamma }{\partial t}$ be a
smooth flow of the partially null curve in $E_{1}^{4\text{ }},$ with the
curvature $k_{3}(s)=0.$ Using definition of $\gamma ,$ we have

\begin{equation}
v^{2}=<\frac{\partial \gamma }{\partial u},\frac{\partial \gamma }{\partial u%
}>.  \tag{3.6}
\end{equation}

Then, by differentiating (3.6), we get%
\begin{equation*}
2v\frac{\partial v}{\partial t}=\frac{\partial }{\partial t}<\frac{\partial
\gamma }{\partial u},\frac{\partial \gamma }{\partial u}>.
\end{equation*}

On the other hand, as $\dfrac{\partial }{\partial u}$ and $\dfrac{\partial }{%
\partial t}$ commute, we have

\begin{equation*}
v\frac{\partial v}{\partial t}=<\frac{\partial \gamma }{\partial u},\frac{%
\partial }{\partial u}\left( \frac{\partial \gamma }{\partial t}\right) >.
\end{equation*}

From (3.3), we obtain

\begin{equation*}
v\frac{\partial v}{\partial t}=<\frac{\partial \gamma }{\partial u},\frac{%
\partial }{\partial u}\left( \beta _{1}T+\beta _{2}N+\beta _{3}B_{1}+\beta
_{4}B_{2}\right) >.
\end{equation*}

By using (2.1), we have%
\begin{equation*}
\frac{\partial v}{\partial t}=<T,\left( \frac{\partial \beta _{1}}{\partial u%
}-k_{1}\beta _{2}v\right) T+\left( \frac{\partial \beta _{2}}{\partial u}%
+\beta _{1}k_{1}v-\beta _{4}k_{2}v\right) N+
\end{equation*}

\begin{equation}
\left( \frac{\partial \beta _{3}}{\partial u}+k_{2}\beta _{2}v\right) B_{1}\
\ +\frac{\partial \beta _{4}}{\partial u}B_{2}>.  \tag{3.7}
\end{equation}

This clearly forces

\begin{equation*}
\frac{\partial v}{\partial t}=\left( \frac{\partial \beta _{1}}{\partial u}%
-\beta _{2}k_{1}v\right) .\ \ \ 
\end{equation*}

\textbf{Lemma 3.4.} Let \ $\{T,N,B_{1},B_{2}\}$ be the Frenet frame of a
partially null curves $\gamma $ and $\dfrac{\partial \gamma }{\partial t}%
=\beta _{1}T+\beta _{2}N+\beta _{3}B_{1}+\beta _{4}B_{2}$ be a smooth flow
of a partially null curves $\gamma $ in $E_{1}^{4\text{ }}.$If $\gamma (s)$
is a unit speed partially null curve in $E_{1}^{4\text{ }},$ with the
curvature $k_{3}(s)=0,$ then we have the following equality

\begin{equation}
\frac{\partial \beta _{1}}{\partial u}=\beta _{2}k_{1}v.\ \ \   \tag{3.8}
\end{equation}

\textbf{Proof. } Let us assume that the partially null curve flow is
inextensible. From (3.4), we have

\begin{equation}
\frac{\partial }{\partial t}s(u,t)=\overset{u}{\underset{0}{\int }}~\frac{%
\partial v}{\partial t}~du=\overset{u}{\underset{0}{\int }}~\left( \frac{%
\partial \beta _{1}}{\partial u}-\beta _{2}k_{1}v\right) ~du=0.\ \  
\tag{3.9}
\end{equation}

This clearly forces%
\begin{equation*}
\frac{\partial \beta _{1}}{\partial u}-\beta _{2}k_{1}v=0.
\end{equation*}

We now restrict ourselves to arc length parametrized curves. That is, $v=1$
and the local coordinate $u$ corresponds to the curve arc length $s$. Then,
we have the following lemma.

\textbf{Lemma 3.5. }Let \ $\{T,N,B_{1},B_{2}\}$ be the Frenet frame of a
partially null curves $\gamma $ and $\frac{\partial \gamma }{\partial t}%
=\beta _{1}T+\beta _{2}N+\beta _{3}B_{1}+\beta _{4}B_{2}$ be a smooth flow
of a partially null curves $\gamma $ in $E_{1}^{4\text{ }}.$ The
differentiations of $\{T,N,B_{1},B_{2}\}$ with respect to $t$ is%
\begin{equation}
\frac{\partial T}{\partial t}=\left( \frac{\partial \beta _{2}}{\partial s}%
+\beta _{1}k_{1}-\beta _{4}k_{2}\right) N+\left( \frac{\partial \beta _{3}}{%
\partial s}+\beta _{2}k_{2}\right) B_{1}+\frac{\partial \beta _{4}}{\partial
s}B_{2},  \tag{3.10}
\end{equation}%
\begin{equation}
\frac{\partial N}{\partial t}=-\left( \frac{\partial \beta _{2}}{\partial s}%
+\beta _{1}k_{1}-\beta _{4}k_{2}\right) T+\psi _{2}B_{1}+\psi _{1}B_{2}, 
\tag{3.11}
\end{equation}%
\begin{equation}
\frac{\partial B_{1}}{\partial t}=-\frac{\partial \beta _{4}}{\partial s}%
T-\psi _{1}N+\psi _{3}B_{1},  \tag{3.12}
\end{equation}%
\begin{equation}
\frac{\partial B_{2}}{\partial t}=-\left( \frac{\partial \beta _{3}}{%
\partial s}+k_{2}\beta _{2}\right) T-\psi _{2}N+\psi _{3}B_{2},  \tag{3.13}
\end{equation}%
where

\begin{equation*}
\psi _{1}=<\frac{\partial N}{\partial t},B_{1}>,\ \ \psi _{2}=<\frac{%
\partial N}{\partial t},B_{2}>,\ \ \ \psi _{3}=<\frac{\partial B_{1}}{%
\partial t},B_{2}>.\ \ 
\end{equation*}

\textbf{Proof. }From the assumption, we have%
\begin{equation*}
\frac{\partial T}{\partial t}=\frac{\partial }{\partial t}\frac{\partial
\gamma }{\partial s}=\frac{\partial }{\partial s}\left( \beta _{1}T+\beta
_{2}N+\beta _{3}B_{1}+\beta _{4}B_{2}\right) .
\end{equation*}

Thus, it is seen that%
\begin{equation*}
\frac{\partial T}{\partial t}=\left( \frac{\partial \beta _{1}}{\partial s}%
+k_{1}\beta _{2}\right) T+\left( \frac{\partial \beta _{2}}{\partial s}%
+\beta _{1}k_{1}-\beta _{4}k_{2}\right) N
\end{equation*}%
\begin{equation}
+\left( \frac{\partial \beta _{3}}{\partial s}+k_{2}\beta _{2}\right) B_{1}+%
\frac{\partial \beta _{4}}{\partial s}B_{2}.  \tag{3.14}
\end{equation}

Subsituting (3.8) into (3.14), we get (3.10).

Since%
\begin{eqnarray*}
&<&T,N>=0~\Rightarrow \ <T,\frac{\partial N}{\partial t}>=-~\left( \frac{%
\partial \beta _{2}}{\partial s}+\beta _{1}k_{1}-\beta _{4}k_{2}\right) , \\
&<&T,B_{1}>=0~\Rightarrow \ <T,\frac{\partial B_{1}}{\partial t}>=-\frac{%
\partial \beta _{4}}{\partial s}, \\
&<&T,B_{2}>=0~\Rightarrow \ <T,\frac{\partial B_{2}}{\partial t}>=-\left( 
\frac{\partial \beta _{1}}{\partial s}+k_{1}\beta _{2}\right) , \\
&<&N,B_{1}>=0~\Rightarrow \ <N,\frac{\partial B_{1}}{\partial t}>=-\psi _{1},
\\
&<&N,B_{2}>=0~\Rightarrow \ <N,\frac{\partial B_{2}}{\partial t}>=-\psi _{2},
\\
&<&B_{1},B_{2}>=1\Rightarrow \ <B_{1},\frac{\partial B_{2}}{\partial t}%
>=-\psi _{3},
\end{eqnarray*}

we have%
\begin{equation*}
<N,\frac{\partial N}{\partial t}>=<B_{1},\frac{\partial B_{1}}{\partial t}%
>=<B_{2},\frac{\partial B_{2}}{\partial t}>=0.
\end{equation*}

In a similar manner above, we can obtain (3.11),(3.12) and (3.13).

\textbf{Theorem 3.6.} Let \ $\{T,N,B_{1},B_{2}\}$ be the Frenet frame of a
partially null curves $\gamma $ and $\dfrac{\partial \gamma }{\partial t}%
=\beta _{1}T+\beta _{2}N+\beta _{3}B_{1}+\beta _{4}B_{2}$ be a smooth flow
of a partially null curves $\gamma $ in $E_{1}^{4\text{ }}.$Then, there
exists the following system of partially differential equation.%
\begin{equation*}
\frac{\partial k_{1}}{\partial t}=\frac{\partial ^{2}\beta _{2}}{\partial
s^{2}}+\frac{\partial \left( \beta _{1}k_{1}\right) }{\partial s}-\frac{%
\partial \left( \beta _{4}k_{2}\right) }{\partial s}-\frac{\partial \beta
_{4}}{\partial s}k_{2}
\end{equation*}

\textbf{Proof \ }From Lemma 3.4, we have%
\begin{eqnarray*}
\frac{\partial }{\partial s}\frac{\partial T}{\partial t} &=&\left( \frac{%
\partial ^{2}\beta _{2}}{\partial s^{2}}+\frac{\partial \left( \beta
_{1}k_{1}\right) }{\partial s}-\frac{\partial \left( \beta _{4}k_{2}\right) 
}{\partial s}\right) N+\left( \frac{\partial \beta _{2}}{\partial s}+\beta
_{1}k_{1}-\beta _{4}k_{2}\right) \left( -k_{1}T+k_{2}B_{1}\right) \\
&&\left( \frac{\partial ^{2}\beta _{3}}{\partial s^{2}}+\frac{\partial
\left( \beta _{2}k_{2}\right) }{\partial s}\right) B_{1}+\frac{\partial
^{2}\beta _{4}}{\partial s^{2}}B_{2}+\frac{\partial \beta _{4}}{\partial s}%
\left( -k_{2}N\right) .
\end{eqnarray*}%
Then\bigskip 
\begin{eqnarray}
\frac{\partial }{\partial s}\frac{\partial T}{\partial t} &=&-\left( \frac{%
\partial \beta _{2}}{\partial s}k_{1}+\beta _{1}k_{1}^{2}-\beta
_{4}k_{1}k_{2}\right) T+\left( \frac{\partial ^{2}\beta _{2}}{\partial s^{2}}%
+\frac{\partial \left( \beta _{1}k_{1}\right) }{\partial s}-\frac{\partial
\left( \beta _{4}k_{2}\right) }{\partial s}-\frac{\partial \beta _{4}}{%
\partial s}k_{2}\right) N  \notag \\
&&+\left( \frac{\partial ^{2}\beta _{3}}{\partial s^{2}}+\frac{\partial
\left( \beta _{2}k_{2}\right) }{\partial s}+\frac{\partial \beta _{2}}{%
\partial s}k_{2}+\beta _{1}k_{1}k_{2}-\beta _{4}k_{2}^{2}\right) B_{1}+\frac{%
\partial ^{2}\beta _{4}}{\partial s^{2}}B_{2}.  \TCItag{3.15}
\end{eqnarray}%
Note that

\begin{equation}
\frac{\partial }{\partial t}\left( \frac{\partial T}{\partial s}\right) =%
\frac{\partial k_{1}}{\partial t}N-\left( \frac{\partial \beta _{2}}{%
\partial s}k_{1}+\beta _{1}k_{1}^{2}-\beta _{4}k_{1}k_{2}\right) T+\psi
_{2}k_{1}B_{1}+\psi _{1}k_{1}B_{2}.  \tag{3.16}
\end{equation}

Hence from (3.15) and (3.16), we get

\begin{equation*}
\frac{\partial k_{1}}{\partial t}=\frac{\partial ^{2}\beta _{2}}{\partial
s^{2}}+\frac{\partial \left( \beta _{1}k_{1}\right) }{\partial s}-\frac{%
\partial \left( \beta _{4}k_{2}\right) }{\partial s}-\frac{\partial \beta
_{4}}{\partial s}k_{2}.
\end{equation*}

This completes the proof.

\textbf{Corollary 3.7.} Let \ $\{T,N,B_{1},B_{2}\}$ be the Frenet frame of a
partially null curves $\gamma $ and $\dfrac{\partial \gamma }{\partial t}%
=\beta _{1}T+\beta _{2}N+\beta _{3}B_{1}+\beta _{4}B_{2}$ be a smooth flow
of a partially null curves $\gamma $ in $E_{1}^{4\text{ }}.$Then,we have the
following equalities.%
\begin{equation*}
\ k_{1}=\frac{1}{\psi _{1}}\left[ \frac{\partial ^{2}\beta _{4}}{\partial
s^{2}}\right] ,
\end{equation*}

\begin{equation*}
\frac{\partial ^{2}\beta _{3}}{\partial s^{2}}+\frac{\partial \left( \beta
_{2}k_{2}\right) }{\partial s}+\frac{\partial \beta _{2}}{\partial s}%
k_{2}+\beta _{1}k_{1}k_{2}-\beta _{4}k_{2}^{2}-\psi _{2}k_{1}=0.
\end{equation*}

\textbf{Theorem 3.8.} Let \ $\{T,N,B_{1},B_{2}\}$ be the Frenet frame of a
partially null curves $\gamma $ and $\dfrac{\partial \gamma }{\partial t}%
=\beta _{1}T+\beta _{2}N+\beta _{3}B_{1}+\beta _{4}B_{2}$ be a smooth flow
of a partially null curves $\gamma $ in $E_{1}^{4\text{ }}.$Then, we have%
\begin{equation}
k_{1}=-\left[ \frac{\partial \psi _{1}}{\partial s}\diagup \frac{\partial
\beta _{4}}{\partial s}\right]   \tag{3.17}
\end{equation}%
and%
\begin{equation}
k_{2}=\frac{1}{\psi _{1}}\left[ \frac{\partial \psi _{3}}{\partial s}\right]
.  \tag{3.18}
\end{equation}

\textbf{Proof. }Noting that $\dfrac{\partial }{\partial s}\left( \dfrac{%
\partial B_{1}}{\partial t}\right) =\dfrac{\partial }{\partial t}\left( 
\dfrac{B_{1}}{\partial s}\right) $, we have $\ $the equations $\left(
3.17\right) $ and $\left( 3.18\right) .$

\textbf{Theorem 3.9.} Let \ $\{T,N,B_{1},B_{2}\}$ be the Frenet frame of a
partially null curves $\gamma $ and $\frac{\partial \gamma }{\partial t}%
=\beta _{1}T+\beta _{2}N+\beta _{3}B_{1}+\beta _{4}B_{2}$ be a smooth flow
of a partially null curves $\gamma $ in $E_{1}^{4\text{ }}.$Then, there
exists the following system of partially differential equation.%
\begin{equation}
\frac{\partial k_{2}}{\partial t}=\frac{\partial \psi _{2}}{\partial s}+%
\frac{\partial \beta _{3}}{\partial s}k_{1}-\beta _{2}k_{1}k_{2}-\psi
_{3}k_{2}.  \tag{3.19}
\end{equation}

\textbf{Proof. \ }By the same\textbf{\ }way above and considering\textbf{\ } 
$\frac{\partial }{\partial s}\left( \frac{\partial B_{2}}{\partial t}\right)
=\frac{\partial }{\partial t}\left( \frac{B_{2}}{\partial s}\right) $ we
reach$.$%
\begin{equation*}
\frac{\partial k_{2}}{\partial t}=\frac{\partial \psi _{2}}{\partial s}+%
\frac{\partial \beta _{3}}{\partial s}k_{1}-\beta _{2}k_{1}k_{2}-\psi
_{3}k_{2}.
\end{equation*}

\bigskip \textbf{4. Inextensible Flows of pseudo null curve in} $E_{1}^{4%
\text{ }}$

We omit the proofs of \ the following theorems became of having close
analogy of the theorems given above.

\textbf{Lemma 4.1.} Let \ $\{T,N,B_{1},B_{2}\}$ be the Frenet frame of a
pseudo null curve $\gamma $ and 
\begin{equation*}
\frac{\partial \gamma }{\partial t}=\alpha _{1}T+\alpha _{2}N+\alpha
_{3}B_{1}+\alpha _{4}B_{2}
\end{equation*}%
be a smooth flow of a pseudo null curve $\gamma $ in $E_{1}^{4\text{ }}.$If $%
\gamma (s)$ be unit speed pseudo null curve in $E_{1}^{4\text{ }},$ with
curvature $k_{1}(s)=1,$ $k_{2}(s)$ and $k_{3}(s)\neq 0,$then we have the
following equality

\begin{equation*}
\frac{\partial v}{\partial t}=\left( \frac{\partial \alpha _{1}}{\partial u}%
-\alpha _{4}v\right) .\ \ \ 
\end{equation*}

\textbf{Lemma 4.2.} Let \ $\{T,N,B_{1},B_{2}\}$ be the Frenet frame of a
pseudo null curve $\gamma $ and $\frac{\partial \gamma }{\partial t}=\alpha
_{1}T+\alpha _{2}N+\alpha _{3}B_{1}+\alpha _{4}B_{2}$be a smooth flow of a
pseudo null curve $\gamma $ in $E_{1}^{4\text{ }}.$If $\gamma (s)$ be unit
speed pseudo null curve in $E_{1}^{4\text{ }},$ with curvature $k_{1}(s)=1,$ 
$k_{2}(s)$ and $k_{3}(s)\neq 0,$then we have the following equality

\begin{equation*}
\frac{\partial \alpha _{1}}{\partial u}=\alpha _{4}v.\ \ \ 
\end{equation*}

\textbf{Lemma 4.3 }Let \ $\{T,N,B_{1},B_{2}\}$ be the Frenet frame of a
pseudo null curve $\gamma $ and $\frac{\partial \gamma }{\partial t}=\alpha
_{1}T+\alpha _{2}N+\alpha _{3}B_{1}+\alpha _{4}B_{2}$ be a smooth flow of a
pseudo null curve $\gamma $ in $E_{1}^{4\text{ }}.$The differentions of $%
\{T,N,B_{1},B_{2}\}$ with respect to $t$ is%
\begin{equation*}
\frac{\partial T}{\partial t}=\left( \frac{\partial \alpha _{2}}{\partial s}%
+\alpha _{1}+\alpha _{3}k_{3}\right) N+\left( \frac{\partial \alpha _{3}}{%
\partial s}+\alpha _{2}k_{2}-\alpha _{4}k_{3}\right) B_{1}+\left( \frac{%
\partial \alpha _{4}}{\partial s}-\alpha _{3}k_{2}\right) B_{2},
\end{equation*}

\begin{equation*}
\frac{\partial N}{\partial t}=-\left( \frac{\partial \alpha _{4}}{\partial s}%
-\alpha _{3}k_{2}\right) T+\psi _{2}N+\psi _{1}B_{1},
\end{equation*}

\begin{equation*}
\frac{\partial B_{1}}{\partial t}=-\left( \frac{\partial \alpha _{3}}{%
\partial s}+\alpha _{2}k_{2}-\alpha _{4}k_{3}\right) T+\psi _{3}N-\psi
_{1}B_{2},
\end{equation*}

\begin{equation*}
\frac{\partial B_{2}}{\partial t}=-\left( \frac{\partial \alpha _{2}}{%
\partial s}+\alpha _{1}+\alpha _{3}k_{3}\right) T-\psi _{3}B_{1}-\psi
_{2}B_{2},
\end{equation*}

where

\begin{equation*}
\psi _{1}=<\frac{\partial N}{\partial t},B_{1}>,\ \ \psi _{2}=<\frac{%
\partial N}{\partial t},B_{2}>,\ \ \ \psi _{3}=<\frac{\partial B_{1}}{%
\partial t},B_{2}>.\ \ 
\end{equation*}

\textbf{Theorem 4.4.} Let \ $\{T,N,B_{1},B_{2}\}$ be the Frenet frame of a
pseudo null curve $\gamma $ and $\frac{\partial \gamma }{\partial t}=\alpha
_{1}T+\alpha _{2}N+\alpha _{3}B_{1}+\alpha _{4}B_{2}$ be a smooth flow of a
pseudo null curve $\gamma $ in $E_{1}^{4\text{ }}.$Then, there exists the
following system of partially differential equations

\begin{equation*}
\frac{\partial ^{2}\alpha _{3}}{\partial s^{2}}+\frac{\partial \left( \alpha
_{2}k_{2}\right) }{\partial s}-\frac{\partial \left( \alpha _{4}k_{3}\right) 
}{\partial s}-\frac{\partial \alpha _{4}}{\partial s}k_{3}+\frac{\partial
\alpha _{2}}{\partial s}k_{2}+\alpha _{1}k_{2}+2\alpha _{3}k_{3}k_{2}-\psi
_{1}=0,
\end{equation*}

\begin{equation*}
\frac{\partial ^{2}\alpha _{2}}{\partial s^{2}}+\frac{\partial \alpha _{1}}{%
\partial s}+\frac{\partial \left( \alpha _{3}k_{3}\right) }{\partial s}+%
\frac{\partial \alpha _{3}}{\partial s}k_{3}+\alpha _{2}k_{2}k_{3}-\alpha
_{4}k_{3}^{2}-\psi _{2}=0,
\end{equation*}

\begin{equation*}
\frac{\partial ^{2}\alpha _{4}}{\partial s^{2}}-\frac{\partial \left( \alpha
_{3}k_{2}\right) }{\partial s}-\frac{\partial \alpha _{3}}{\partial s}%
k_{2}-\alpha _{2}k_{2}+\alpha _{4}k_{2}k_{3}=0,
\end{equation*}%
\begin{equation*}
\frac{\partial k_{3}}{\partial s}=\frac{\partial \psi _{3}}{\partial s}-%
\frac{\partial \alpha _{3}}{\partial s}-\alpha _{2}k_{2}+\alpha
_{4}k_{3}-\psi _{2}k_{3},
\end{equation*}%
\begin{equation*}
\frac{\partial k_{2}}{\partial s}=\frac{\partial \psi _{1}}{\partial s}+\psi
_{2}k_{2},
\end{equation*}

\begin{equation*}
\frac{\partial \psi _{2}}{\partial s}=\frac{\partial \alpha _{4}}{\partial s}%
-\alpha _{3}k_{2}-\psi _{1}k_{3}+\psi _{2}k_{3}.
\end{equation*}

\textbf{References}

[1] \ H.Balgetir, M.Bekta\c{s}, M. Erg\"{u}t, On a characterization of null
helices, Bull. Inst. Math. Acad. Sinica 29, 71-78, 2001.

[2] \ H.Balgetir, M. Bekta\c{s} and J.Inoguchi, Null Bertrand curves and
their characterizations, Note Mat. 23, no. 1, 7-13, .2004

[3] \ Barros, M., General helices and a theorem of Lancret, Proc. Amer.
Math. Soc. 125,1503-1509, 1997.

[4\ ]\ S. Bas and T. K\"{o}rp\i nar, Inextensible Flows of Spacelike Curves
on Spacelike Surfaces according to Darboux Frame in M31 , Bol. Soc. Paran.
Mat. 31 (2) , 9--17, 2013.

[5] \ S. Bas , T. K\"{o}rp\i nar and E. Turhan, New Type Inextensible Flows
of Timelike Curves \.{I}n Minkowski space-time $M_{1}^{4},$ AMO Advanced
Modeling and Optimization, Volume 14, Number 2, 2012.

[6] \ G. Chirikjian, J. Burdick, A modal approach to hyper-redundant
manipulator kinematics,IEEE Trans. Robot. Autom. 10, 343--354,1994.

[7] \ M. Desbrun, M.P. Cani-Gascuel, Active implicit surface for animation,
in: Proc. Graphics Interface Canadian Inf. Process. Soc., 143--150, 1998

[8] \ A.Ferrandez, A. Gimenez and P. Lucas, Null helices in Lorentzian space
forms, International Journal of Modern Physics A 16, 4845--4863, 2001.

[9] \ N. G\"{u}rb\"{u}z, Inextensible flows of spacelike, timelike and null
curves, Int. J. Contemp. Math. Sciences, Vol. 4, no. 32, 1599-1604, 2009.

[10] \ K.Ilarslan and E.Nesovic, Some Characterizations of Null, Pseudo Null
and Partially Null Rectifiying Curves in Minkowski\ Space-Time, Taiwanese
Journal of Mathematics,Vol 12, No:5,pp 1035-1044, 2008.

[11] \ M. Kass, A. Witkin, D. Terzopoulos, Snakes: active contour models,
in: Proc. 1st Int. Conference on Computer Vision, 259--268,1987.

[12]\ \ \ T.K\"{o}rp\i nar and E.Turhan, New Inextensible Flows of Timelike
Curves on the Oriented Timelike Surfaces According to Darboux Frame \.{I}n M$%
_{1}^{3},$ AMO \TEXTsymbol{\vert} Advanced Modeling and Optimization, Volume
14, Number 2, 2012.

[13]\ \ \ T.K\"{o}rp\i nar and E.Turhan,\ A New Version of Inextensible
Flows of Spacelike Curves with Timelike B2 in Minkowski Space-Time E41
,Differ Equ Dyn Syst, DOI 10.1007/s12591-012-0152-4

[14] M.K\"{u}lahc\i , M.Bekta\c{s}, M.Erg\"{u}t, Curves of AW(k)-type in
3-dimensional Null Cone, Physics Letters A. 371, 275-277, 2007.

[15] M.K\"{u}lahc\i , M.Bekta\c{s},M.Erg\"{u}t, On Harmonic Curvatures of
Null Curves of the AW(k)-Type in Lorentzian Space Z. Naturforsch. 63a, 248
-- 252, 2008.

[16\ ] \ D. Y. Kwon, F.C. Park, D.P. Chi, Inextensible flows of curves and
developable surfaces, Appl. Math. Lett. 18, 1156-1162, 2005.

[17] D.Latifi, \ A.Razavi A , Inextensible Flows of Curves in Minkowskian
Space, Adv. Studies Theor. Phys. 2(16): 761-768, 2008.

[18] \ H.Q. Lu, J.S. Todhunter, T.W. Sze, Congruence conditions for
nonplanar developable surfaces and their application to surface recognition,
CVGIP, Image Underst. 56, 265--285, 1993.

[19] \ H. Mochiyama, E. Shimemura, H. Kobayashi, Shape control of
manipulators with hyper degrees of freedom, Int. J. Robot.Res., 18,
584--600, 1999.

[20] \ A.O. Ogrenmis, M. Yenero\u{g}lu, Inextensible curves in the Galilean
Space, International Journal of the Physical Sciences, 5(9),1424-1427, 2010.

[21] D.J. Unger, Developable surfaces in elastoplastic fracture mechanics,
Int. J. Fract. 50, 33--38, 1991.

[22] \"{O} . G. Y\i ld\i z, S. Ersoy, M. Masal, A note on inextensible flows
of curves on oriented surface, arXiv:1106.2012v1.

[23] \ \"{O} . G. Y\i ld\i z, M. Tosun, S. O. Karaku\c{s}, A note on
inextensible flows of curves in $E^{n}$, arXiv:1207.1543v1

[24] \"{O}.G. Y\i ld\i z and M.Tosun, A Note on Inextensible Flows of Curves
in $E_{1}^{n},$ arXiv:1302.6082v1 [math.DG] 25 Feb 2013

\bigskip

\bigskip

\bigskip

\end{document}